\documentclass[a4paper,11pt]{article}

\usepackage{url}
\usepackage{blkarray}
\usepackage{hyperref,pgf}
\usepackage{setspace}
\usepackage[version=3]{mhchem}
\usepackage{amssymb,amsfonts,amsmath,amsthm}
\usepackage{xypic,enumerate}
\usepackage{color,subfig,tikz,pgf}
\usetikzlibrary{snakes}
\usetikzlibrary{arrows} 
\usepackage{cite} 
\input xy
\xyoption{all}

\newcommand{\R}{\mathbb{R}}
\newcommand{\Z}{\mathbb{Z}}

\renewcommand{\k}{\kappa}

\DeclareMathOperator{\im}{im} 

\DeclareMathOperator{\diag}{diag}
\DeclareMathOperator{\rank}{rank}

\newcommand{\mmF}{\mathcal{F}}
\newcommand{\mmG}{\mathcal{G}}

\newcommand{\st}{\mid} 
\renewcommand{\k}{\kappa}
\newcommand{\wk}{\widehat{\kappa}}

\newenvironment{enumerate*}[1][{}]{\begin{itemize}\renewcommand{\labelitemi}{#1}}{\end{itemize}}

\newtheorem{proposition}{\bf Proposition}[section]

\newtheorem{definition}{\bf Definition}[section]
\newenvironment{defi}[1][{}]{\begin{definition}\normalfont{#1}}{\end{definition}}

\makeatletter
\def\blfootnote{\xdef\@thefnmark{}\@footnotetext}
\makeatother

\begin{document}

\title{Injectivity, multiple zeros, and multistationarity in reaction networks}

\author{%
Elisenda Feliu
}

\blfootnote{
\scriptsize
\noindent
{\bf E. Feliu} (\href{efeliu@math.ku.dk}{efeliu@math.ku.dk}). 
Dept.\ of Mathematical Sciences, University of Copenhagen, Universitetsparken 5, 2100 Denmark 
}

\maketitle

\begin{abstract}
Polynomial dynamical systems are widely used to model and study real phenomena. In biochemistry, they are the preferred choice for modelling the concentration of chemical species in reaction networks with mass-action kinetics.
These systems are typically parameterised by many (unknown) parameters. A  goal is to understand how properties of the dynamical systems depend on the parameters.

Qualitative properties relating to the behaviour of a dynamical system are locally inferred from the system at steady state. 
Here we focus on steady states that  are the positive solutions to a parameterised system of \emph{generalised polynomial} equations. In recent years, methods from computational algebra have been developed to understand these solutions, but our knowledge is limited: for example, we cannot efficiently decide how many positive solutions the system has as a function of the parameters. Even deciding whether there is one or more solutions is non-trivial.

We present a new method, based on  so-called \emph{injectivity}, to preclude or assert that multiple positive solutions exist. The results apply to generalised polynomials and variables can be restricted to the  linear, parameter-independent  first integrals of the  dynamical system. 
The method has been  tested in a  wide range of systems.
\end{abstract}

\section{Introduction}
The cell's ability to respond on-off to gradual changes in a signal  is linked to bistability,  the ability of a dynamical system to admit two different stable steady states  \cite{Markevich-mapk,TG-nature}. Specifically, bistability underlies the emergence of hysteresis and switch-like behavior.

In the context of molecular biology, models usually depend on unknown parameters and one is interested in determining whether there exists a choice of parameter values for which the system has two stable steady states.
In particular, the concentration of the species in a \emph{(bio)chemical reaction network} is typically modelled with 
a system of ordinary differential equations (ODEs):
$$\frac{dx}{dt} = f_\k(x),\qquad x\in \R^n_{\geq 0},$$
where $\k$ is a vector of parameters. Bistability arises when the equation $f_\k(x)=0$ admits two (positive) solutions that correspond to stable steady states. For unknown parameters, determining whether a steady state is stable is highly non-trivial. Therefore, the focus has mainly been on determining whether the equation $f_\k(x)=0$ admits multiple solutions for some parameter values $\k$.

Determining the  number and stability of steady states has been an important problem in the context of chemical reaction networks in the 20th century, and there is a rich literature covering both theoretical and experimental aspects. The majority of the studies date back to the seventies and eighties (see for example the review \cite{razon:1987dv} and the references therein).   
In particular, several foundational theoretical studies on determining the number and stability of steady states of generic chemical reaction networks were developed in that period. These studies address different broad classes of reaction networks, e.g. detailed or complex balanced networks or reaction networks  with a specific form. 
Of these foundational works we highlight the so-called Chemical Reaction Network Theory (CRNT) of Jackson, Horn, and Feinberg, e.g. \cite{feinberg-def0,Feinberg1988,feinberg-invariant,crnttoolbox,horn,hornjackson}, the work of Clarke, e.g. \cite{Clarke:1980tz,Clarke:1975jl}, the  work of Vol'pert, e.g. \cite{volpert,volpert2} and several other works \cite{GORBAN:1986kia,Fedotov:1984hg,ivanova,bykov2}.

Since the advent of systems biology in the late nineties, complex and high-dimensional reaction networks involving biochemical entities have become the object of study. In this context  the methods developed for chemical reaction networks are often not applicable, due to the complexity of the networks. For instance, most networks modelling cell signalling  events are neither detailed nor complex balanced. A new generation of methods has been developed, often within the framework of CRNT, incorporating new formalism and approaches from computational algebra, e.g. \cite{bykov,PerezMillan,wiuf-feliu,gatermannbook, conradi-PNAS}.

The extensive list of existing methods apply under different assumptions and to different modelling strategies. Most methods  specialise to the steady states  with positive  coordinates, that is, 
 the methods decide whether the function $f_\k$ admits multiple positive zeros and exclude the zeros at the boundary. Unfortunately, there does not exist a  single method that can satisfactorily handle any given reaction network. 
 
In this paper, we develop a new method to address the existence or preclusion of multiple steady states. The method builds on so-called injectivity-based methods, which underly several earlier approaches, for instance \cite{craciun-feinbergI,banaji-craciun2,soule,muller,PerezMillan,wiuf-feliu,shinar-conc}. For polynomial equations, the setting  was introduced in \cite{craciun-feinbergI}, and subsequently extended in \cite{craciun-feinberg-semiopen,feliu-inj,wiuf-feliu,MullerSigns,CGS}.

\medskip
The idea underlying injectivity methods is simple: Consider a map
$f\colon \R^n_{>0}\rightarrow \R^n$. If $f$ is  injective, then there do not exist distinct $x,y$ in $\R^n_{>0}$ such that $f(x)=f(y)$.
In particular, $f$ does not have two  distinct zeros in $\R^n_{>0}$.

The converse is  false. For instance, consider the map $f\colon \R^2_{>0}\rightarrow \R^2$
given by 
\begin{align*}
f_1(x_1,x_2) &=x_1x_2-10 x_1+0.1x_2,\\
 f_2(x_1,x_2) & =x_1x_2+0.1x_1-10x_2. 
 \end{align*}
The only positive zero of $f$ is $(9.9,9.9)$, but $f$ is not injective on $\R^2_{>0}$ since
$$f(5.9,2)=f(7.9,4)=(-47,-7.61). $$
Thus, in this case, we cannot conclude that $f$ does not have multiple positive zeros by checking whether $f$ is injective.
However, we can make another simple observation: The positive zeros of $f$ are also the positive solutions to:
\begin{equation}\label{eq:g1}
x_1-10 x_1x_2^{-1}= -0.1 ,\quad x_2-10x_1^{-1}x_2=-0.1. 
\end{equation}
An easy check shows that the map $g\colon \R^2_{>0}\rightarrow \R^2$ defined by
\begin{equation}\label{eq:g3}
g(x_1,x_2) =(x_1-10 x_1x_2^{-1},x_2-10x_1^{-1}x_2)
 \end{equation}
 is injective on $\R^2_{>0}$ and hence we conclude that neither \eqref{eq:g1} nor $f=0$ have multiple positive solutions. Observe that positivity of the solutions is critical here, because $g$ is obtained by dividing its two components by $x_2$ and $x_1$, respectively.

\medskip
This simple example illustrates the approach introduced here: By checking the injectivity of a new map $g$ constructed from $f$, we  conclude that $f$ does not have more than one positive zero. 
We focus on the positive zeros of so-called \emph{generalised polynomial maps},  that is, polynomials with real exponents. For example, $g$ in \eqref{eq:g3} is a generalised polynomial map. 
More generally, we consider families of parameterised generalised polynomial maps (cf. \eqref{eq:FAV}) and use the idea above to preclude or determine the existence of multiple zeros on the positive part of affine subspaces $x^*+S$, where $x^*\in \R^n_{>0}$ and  $S\subseteq \R^n$ is a vector subspace.

The method is designed to address  \emph{multistationarity}  in reaction networks (see Section \ref{sec:net}). However, the method  focuses generally on the positive zeros of generalised polynomial maps. Therefore, we start by  presenting the method in a mathematical framework in Section \ref{sec:math}. In particular we  explain how injectivity can be  used both to preclude and to assert the existence of multiple positive zeros. We proceed to exemplify how the method can be used to address multistationarity for a wide range of reaction networks in Section~\ref{sec:net}. The steps of the method are summarised in subsection~\ref{sec:math}\ref{sec:method}.

\medskip
Given a positive integer $n$, we let $[n]:=\{1,\dots,n\}$.

\section{Mathematical framework}\label{sec:math}

\subsection{Injectivity and multiple zeros}

\paragraph{S-injectivity. }
Given a vector subspace $S\subseteq \R^n$, we first define the notion of injectivity on the positive part of the affine subspaces $x^*+S$ with $x^*\in \R^n_{>0}$, cf. \cite{MullerSigns}.

\begin{defi}\label{def:Sinj}
Let $f\colon \R^n_{>0}\rightarrow \R^m$ be a map and $S\subseteq \R^n$ a vector subspace. We say that
\begin{itemize}
\item $f$ is \emph{$S$-injective} if for all distinct $x,y\in \R^n_{>0}$ such that $x-y\in S$ we have $f(x)\neq f(y)$.  
\item $f$ has \emph{multiple $S$-zeros} if there exist distinct $x,y\in \R^n_{>0}$ such that $x-y\in S$ and $f(x)=f(y)=0$.
\end{itemize}
\end{defi}
Equivalently, $f$ does not have multiple $S$-zeros  if
$f$ does not have two  distinct zeros  in some set $(x^*+S)\cap \R^n_{>0}$ with $x^*\in \R^n_{>0}$.
Clearly, if $f$ is $S$-injective, then $f$ does not have multiple $S$-zeros.

\begin{defi}
Let $S\subseteq \R^n$ be a vector subspace and $\mmF \subseteq \{f\colon \R^n_{>0} \rightarrow \R^m \}$ be a set of functions. We say that
\begin{itemize}
\item   $\mmF$ is $S$-injective, if $f$ is $S$-injective for all $f \in \mmF$.
\item  $\mmF$ has multiple $S$-zeros, if there exists  $f\in \mmF$ such that $f$ has multiple $S$-zeros.
\end{itemize}
\end{defi}

\paragraph{Generalised polynomial maps. }
We consider families of generalised polynomial maps given by two real matrices:   $A=(a_{ij})\in \R^{m\times r}$ (coefficient matrix) and  $V=(v_{ij})\in \R^{n\times r}$ (exponent matrix). We denote the $j$-th column of $V$ by $v^j$ and define, for  $x\in \R^n_{>0}$,  the vector $x^V\in \R^r_{>0}$ as
$$\big(x^V\big)_j = x^{v^j} = x_{1}^{v_{1j}}\cdot \ldots \cdot x_{n}^{v_{nj}},\quad j\in [r].$$
The vector $x^V$ is a vector of generalised monomials, whose exponents are given by the columns of $V$. For a positive vector $\kappa\in \R^r_{>0}$, we define the map $f_\kappa\colon \R^n_{>0}\rightarrow \R^m$ as 
\begin{equation}\label{eq:f}
 f_\kappa(x) = A\big(\kappa \circ x^V\big),
 \end{equation}
 where $\circ$ denotes the  product componentwise (or Hadamard product). The coefficient $\k_j$ of each monomial  will be treated as a parameter. 
Since 
$$f_\kappa(x)=f_\kappa(y)\quad \textrm{ if and only if }\quad A\big(\kappa \circ \big(x^V-y^V\big)\big)=0,$$ then $S$-injectivity of  $f_\kappa$ 
is unchanged if the matrix $A$ is replaced by a matrix $\widetilde{A}$ with the same kernel as $\widetilde{A}$.
 
 \medskip
For example, the map
\begin{equation}\label{eq:ex1}
f_\k(x) = \begin{pmatrix}  \k_1 x_1 x_2 - \k_2 x_2x_3 + \k_3 x_1 x_4 \\   \k_2 x_2x_3 - \k_3 x_1 x_4 - \k_4 x_2 x_4  \end{pmatrix},
\end{equation}
 takes the form \eqref{eq:f} with 
\begin{equation}\label{eq:AVex1}
A =  {\small \left(\begin{array}{rrrr}  1 & -1 & 1 & 0  \\  0 & 1 &  -1 & -1  \end{array}\right)},  \qquad 
V={\small \left(\begin{array}{rccc} 1 & 0 & 1 & 0  \\  1 & 1 & 0 & 1   \\  0 & 1 & 0 & 0  \\   0 & 0 & 1 & 1 \end{array}\right)},  
\end{equation}
such that
$\kappa \circ x^V = (\k_1 x_1x_2,\k_2 x_2x_3,\k_3 x_1x_4,\k_4 x_2x_4)^t. $

\medskip
Any generalised polynomial map can be written in the form \eqref{eq:f}, possibly for more than one choice of   $A$, $V$, and  $\k$.
For instance,  a generalised polynomial map $f\colon \R^n_{>0} \rightarrow \R^m$ involving monomials $x^{v^1},\dots,x^{v^r}$  can be written componentwise  as
$$f_i(x)= \sum_{j=1}^r a_{ij} x^{v^j},\qquad a_{ij}\in \R, \quad i\in [n]. $$
Hence, $f$ is of the form \eqref{eq:f} with $\kappa=(1,\dots,1)^t$, $A=(a_{ij})$ and $V$ having columns $v^1,\dots,v^r$.

We let
\begin{equation}\label{eq:FAV}
\mmF_{A,V} := \big\{ f_\k\colon \R^n_{>0} \rightarrow \R^m\st f_\k(x) =A \big(\kappa \circ x^V\big),  \kappa\in \R^r_{>0} \big\} 
\end{equation}
be the family of all generalised polynomial maps \eqref{eq:f} obtained by varying the vector of parameters $\kappa$.
Here, we address the question of  determining whether $\mmF_{A,V}$ has multiple $S$-zeros.
If the family $\mmF_{A,V}$ is $S$-injective, then $\mmF_{A,V}$ does not have multiple $S$-zeros, but the converse is not true.

In \cite{MullerSigns}, a characterisation of the families of generalised polynomial maps $\mmF_{A,V}$ that are $S$-injective is given in terms of sign vectors and symbolic determinants. 
Namely, the family is $S$-injective if and only if some particular sets of signs do not intersect. These sets are constructed from the orthants of $\R^n$ and $\R^r$ that certain subspaces, defined from $A,V,S$, intersect.

When $\dim(S)\geq \rank(A)$, $S$-injectivity of $\mmF_{A,V}$ can  also be decided based on the signs of the coefficients of the determinant of a symbolic matrix. The latter scenario is applicable when studying steady states of chemical reaction networks. Due to its simplicity and applicability, we describe only the determinant-based criterion for injectivity here.  The reader is referred to 
\cite{MullerSigns} for the sign criterion.

\paragraph{Determinant criterion for injectivity}
Let $S\subseteq \R^n$ be a vector subspace of dimension $s$ and $A\in \R^{m\times r}$ such that $s\geq \rank(A)$. 
 Let $Z\in \R^{(n-s)\times n}$ be any matrix whose rows  are a basis of $S^\perp$ and  choose  $\widetilde{A}\in \R^{s\times r}$ such that $\ker(A)=\ker\hspace{-0.05cm}\big(\widetilde{A}\big)$. For example, we can choose a set of $s$ rows of $A$ with the same rank as $A$.
 
 For $\k\in \R^r_{>0}$ and $\lambda\in \R^n_{>0}$, let 
 $M_{\k,\lambda}$ be the square matrix given in block form as
 \begin{equation}\label{eq:mkl}
 M_{\k,\lambda} = \begin{pmatrix}  Z \\ \widetilde{A} \diag(\k) V^t \diag(\lambda) \end{pmatrix},
 \end{equation}
  where, for a vector $w$,  $\diag(w)$ denotes the diagonal matrix with diagonal $w$.
    By considering $\k,\lambda$ as indeterminates, the determinant of $M_{\k,\lambda}$ is a polynomial $p(\k,\lambda)$ in $\k,\lambda$. 
    It is a result of \cite{wiuf-feliu,MullerSigns}, that the family $\mmF_{A,V}$ is $S$-injective if and only if $p(\k,\lambda)$ is not identically zero and, further, all its coefficients have the same sign.

Note than when $\rank(A)<s$, then $p(\k,\lambda)$ is identically zero and the family $\mmF_{A,V}$ is not $S$-injective.

\subsection{Gauss reduction and injectivity}\label{sec:inj}

In this subsection we present the steps of the method that can be used to conclude that $\mmF_{A,V}$ does not have multiple $S$-zeros.

\paragraph{The family $\mmG_{A,V}$. }
Let $S\subseteq \R^n$ be a vector subspace of dimension $s$ and  assume that $A$ has rank $s$. 
We choose as above a matrix $\widetilde{A}\in \R^{s\times r}$ such that $\ker(A)=\ker\big(\widetilde{A}\big)$ and let $B$ be the Gauss reduction of $\widetilde{A}$.  
Since the kernel of $A$ and $B$ agree, the family $\mmF_{B,V}$ is $S$-injective or has multiple $S$-zeros  if and only if this is the case for $\mmF_{A,V}$. 
By reordering the columns if necessary, we can assume that  columns  $1$ to $s$ of $\widetilde{A}$ are linearly independent. 
Therefore, without loss of generality, we restrict   to  families $\mmF_{A,V}$ with 
\begin{equation}\label{eq:B}
A=  \begin{pmatrix} id_s\ | \ A_1\end{pmatrix}, 
\end{equation}
where $id_s$ is the identity matrix of size $s$ and $A_1\in \R^{s\times (r-s)}$.
Then $f_\k\in \mmF_{A,V}$ satisfies $f_\k(x)=0$ if and only if
$$  \begin{pmatrix} \k_1 x^{v^1} \\ \vdots \\ \k_s x^{v^s} \end{pmatrix}  = - A_1\hspace{-0.1cm}   \begin{pmatrix} \k_{s+1} x^{v^{s+1}} \\ \vdots \\ \k_r x^{v^r}   \end{pmatrix},$$ 
and this equality holds for $x\in \R^n_{>0}$ if and only if 
\begin{equation}\label{eq:geta}
  \begin{pmatrix} \k_1 \\ \vdots \\ \k_s  \end{pmatrix} = g_{\wk}(x),\quad \textrm{where } \quad g_{\wk}(x)= - \begin{pmatrix} x^{-v^1} \\  \vdots \\ x^{-v^s}\end{pmatrix} \circ A_1 \hspace{-0.1cm}  \begin{pmatrix} \k_{s+1} x^{v^{s+1}} \\ \vdots \\ \k_r x^{v^r}   \end{pmatrix},
 \end{equation}
with $\wk=(\k_{s+1},\dots,\k_r)$. 
If the map $g_{\wk}(x)$ is $S$-injective, then $f_\k$ does not have multiple $S$-zeros. 
Therefore, by defining
$$\mathcal{G}_{A,V} := \{ g_{\wk}\colon \R^n_{>0} \rightarrow \R^s \st  \wk \in \R^{r-s}_{>0} \}, $$
we have that $\mmF_{A,V}$ does not have multiple $S$-zeros if $\mathcal{G}_{A,V}$ is $S$-injective.

\paragraph{Example 1. } Consider the family $\mmF_{A,V}$ defined by \eqref{eq:AVex1} and 
let 
\begin{equation}\label{eq:S}
S=\langle(1,-1,0,0),(0,0,1,-1)\rangle\subseteq \R^4.
\end{equation} 
Both the dimension of $S$ and the rank of $A$ are $2$. Hence we take $\widetilde{A}=A$ and replace $A$ by its Gauss reduction:
\begin{equation}\label{eq:Ared}
A=\left(\begin{array}{rrrr}  1 & 0 & 0 & -1  \\  0 & 1 &  -1 & -1  \end{array}\right).
\end{equation}
The determinant criterion tells us that the family $\mmF_{A,V}$ is not $S$-injective.
The map $g_{\wk}$ in \eqref{eq:geta} is
\begin{align}
g_{\k_3,\k_4}(x) &= - \begin{pmatrix} x_1^{-1}x_2^{-1}\\   x_2^{-1}x_3^{-1} \end{pmatrix} \circ \left(\begin{array}{rr}   0 & -1  \\   -1 & -1  \end{array}\right)  \begin{pmatrix}  \k_3 x_1 x_4  \\  \k_4 x_2 x_4 \end{pmatrix}  = \begin{pmatrix}  \k_4 x_1^{-1}x_4  \\   \k_3 x_1 x_2^{-1}x_3^{-1} x_4 + \k_4  x_3^{-1}x_4 \end{pmatrix} \nonumber \\ 
&=  \begin{pmatrix}  1 & 0 & 0  \\  0 & 1 & 1 \end{pmatrix}   \begin{pmatrix}  \k_4 x_1^{-1} x_4  \\  \k_3 x_1 x_2^{-1}x_3^{-1} x_4  \\  \k_4  x_3^{-1}x_4 \end{pmatrix}. \label{eq:geta2}
\end{align}
The family $\mmG_{A,V} $ is not of the form \eqref{eq:FAV}, since two different monomials have the same parameter as coefficient.
However, observe that 
$$\mmG_{A,V} = \left\{ f_\eta \in \mmF_{A',V'} \st \eta\in  \R_{>0}^3,\   \eta_1=\eta_3 \right\}\subseteq \mmF_{A',V'}$$
with 
\begin{equation}\label{eq:CW}
A'={\small \begin{pmatrix}  1 & 0 & 0  \\  0 & 1 & 1 \end{pmatrix}},\qquad V'= {\small  \left(\begin{array}{rrr}  -1 & 1 & 0   \\  0 & -1 & 0  \\ 0 &  -1 & -1 \\ 1 & 1 & 1  \end{array} \right)}.  
\end{equation}
We use the determinant criterion to decide whether the family $\mmF_{A',V'}$ is $S$-injective. For suitable $Z$, the matrix $M_{\k,\lambda}$ is 
$$M_{\k,\lambda} = {\small  \begin{pmatrix}  1 & 1 & 0 & 0  \\  0 & 0 & 1 & 1 \\ -\k_1\lambda_1 & 0 & 0 & \k_1\lambda_4 \\
\k_2\lambda_1 & -\k_2\lambda_2 & -(\k_2+\k_3)\lambda_3 & (\k_2+\k_3)\lambda_4
\end{pmatrix}}.  $$ 
The determinant of this matrix is 
$$\det(M_{\k,\lambda} )=-\k_{1}(\k_{2}\lambda_{1}\lambda_{3}+2\,\k_{2}\lambda_{1}\lambda_{4}+\k_{2}\lambda_{2}\lambda_{4}+\k_{3}\lambda_{1}\lambda_{3}+\k_{3}\lambda_{1}\lambda_{4}).
$$
Since all coefficients have the same sign, the family $\mmF_{A',V'}$   is $S$-injective.
Hence, so is the family $\mmG_{A,V}$ and it follows that  $\mmF_{A,V}$ does not have multiple $S$-zeros.

\paragraph{The family $\mmF_{A',V'}$. }
Although the family $\mathcal{G}_{A,V}$ might not be of the type \eqref{eq:FAV},  one can always find a family  $\mmF_{A',V'}$, with $A'\in \R^{s\times r'}$,  $V'\in \R^{n\times r'}$, and such that $$\mmG_{A,V}\subseteq \mmF_{A',V'}.$$ 
Indeed, one first identifies the different generalised monomials in $\wk$, and $x$ of $g_{\wk}(x)$.
The exponents of the monomials in $x$ define $V'$. The matrix $A'$ is the matrix of coefficients. When doing so, some parameter $\k_j$ might be multiplying two different generalised monomials in $x$, say $x^{w_1}$ and $x^{w_2}$, as it is the case in the example. 
This occurs whenever a column $j$ of $A_1$ contains  two nonzero entries at rows $i_1,i_2$ such that $v^{i_1}\neq v^{i_2}$.
Then the monomial $\k_{s+j}x^{v^{s+j}}$ appears as a summand in two different rows of the product 
$$ A_1\big(\kappa \circ \big(x^{v^{s+1}}, \ldots, x^{v^r}\big)\big)^t,$$
and these rows are multiplied by different monomials $x^{-v^{i_1}},x^{-v^{i_2}}$ when constructing $g_{\wk}$.
In general, there exists a partition $[r']= I_1\cup \dots \cup I_q$ and 
matrices $A'\in \R^{s\times r'}$ and $V'\in \R^{n\times r'}$ such that  
\begin{equation}\label{eq:I}
\mmG_{A,V} = \big\{f_\eta \in  \mmF_{A',V'} \st  \eta\in \R^{r'}_{>0},\  \eta_i=\eta_j \textrm{ if }i,j\in I_k, \textrm{ for some }k \in [q]\big\}.
\end{equation}
In the example above, $r'=3$ and we let $I_1=\{2\}$, $I_2=\{1,3\}$.

\paragraph{Summary of the steps. }
Given a family of generalised polynomial maps $\mmF_{A,V}$ such that $A$ is Gauss reduced and of the form \eqref{eq:B}, then we construct the families $\mmG_{A,V}$ and $\mmF_{A',V'}$.
Since the matrix  $A'$ has $s$ rows, we apply the determinant criterion to determine whether  $ \mmF_{A',V'}$ is  $S$-injective.
If this is the case, then we conclude that $\mmF_{A,V}$ does not have multiple $S$-zeros.
 
Observe that when constructing  $\mmG_{A,V}$, the parameters $\k_1,\dots,\k_s$ are selected and ``eliminated''. 
We could have selected any other set of $s$ parameters $\k_{i_1},\dots,\k_{i_s}$, as long as the columns ${i_1},\dots,i_s$ of $A$ are linearly independent.  In this case, the  parameters $\k_{i_1},\dots,\k_{i_s}$ are eliminated if we  first permute
the columns of $A$ and $V$, and the entries of $\k$ simultaneously such that the indices ${i_1},\dots,i_s$ are the first indices. We then apply Gauss reduction and proceed as above  with the new data.
 
 Therefore, if $\mmF_{A',V'}$ is not $S$-injective, then we repeat the described process after simultaneously permuting the columns of $A$ and $V$, and the entries of $\k$. We do this for all possible permutations such that the first $s$ columns of $A$ are linearly independent, unless for some permutation we conclude that $\mmF_{A,V}$ does not have multiple $S$-zeros and the procedure stops.

\paragraph{The family $\mmF_{\widehat{A},\widehat{V}}$. }
In preparation for the next section, in particular for the proof of Proposition~\ref{prop:g}, we embed $\mmG_{A,V}$ into an even larger family $\mmF_{\widehat{A},\widehat{V}}$ where:
\begin{itemize}
\item The matrix
$\widehat{A}\in \R^{s\times s(r-s)}$ is the block diagonal matrix with diagonal blocks given by minus the rows of $A_1=(\overline{a}_{ij})$:   $$\widehat{a}_{i,(i-1)(r-s)+j}= - \overline{a}_{ij},\qquad\textrm{for}\quad i\in [s],\ j\in [r-s],$$ 
and zero otherwise. 
\item 
The $\ell$-th column of  $\widehat{V}\in \R^{n\times s(r-s)}$ is 
$$v^{j+s} - v^i,\qquad\textrm{if}\quad\ell=(i-1)(r-s)+j,\quad\textrm{with}\quad i\in [s], \ j\in [r-s].$$
\end{itemize}

Given $\varphi=f_\eta \in \mmF_{\widehat{A},\widehat{V}}$ with $\eta\in\R^{s(r-s)}_{>0}$, then $\varphi \in \mmG_{A,V}  $
if and only if 
$$\eta_{j+(b-1)(r-s)} = \eta_{j+(c-1)(r-s)}$$
for all $j\in [r-s]$ and $b,c\in [s]$ such that $\overline{a}_{bj}\overline{a}_{cj}\neq 0$.

For each $j\in [r-s]$, let $\alpha_j\subseteq \{1,\dots,s\}$ be the support of the $j$-th column of $A_1$.
Then,  
\begin{equation}\label{eq:AV}
\mmG_{A,V} = \{ f_\eta \in \mmF_{\widehat{A},\widehat{V}} \st \eta_{j+(b-1)(r-s)} = \eta_{j+(c-1)(r-s)}, \textrm{ if }b,c\in \alpha_j  \}.
\end{equation}
Observe that only the sets $\alpha_j$ with cardinality at least $2$ are relevant.
Further, we have
$$\mmG_{A,V}  \subseteq  \mmF_{A',V'} \subseteq \mmF_{\widehat{A},\widehat{V}}.$$

\medskip
For example, consider the matrix $V$ in \eqref{eq:AVex1} and $A$ in \eqref{eq:Ared}, where $s=2$, $r=4$, and $n=4$. We have
$$
- A_1=\left(\begin{array}{rrrr}   0 & 1  \\     1 & 1  \end{array}\right),\qquad \widehat{A}=\left(\begin{array}{rrrr}   0 & 1 &0&0 \\    0&0& 1 & 1  \end{array}\right),\qquad \textrm{and}\qquad 
\widehat{V} ={\small \left(\begin{array}{rrrr} 0 & -1 & 1 & 0  \\  -1 & 0 & -1 & 0   \\  0 & 0 & -1 & -1  \\   1 & 1 & 1 & 1 \end{array}\right)}.$$
A map $f_\eta\in \mmF_{\widetilde{A},\widetilde{V}}$ is of the form
\begin{align*}
 f_\eta(x) &= \widetilde{A}\big( (\eta_1x_2^{-1}x_4,\eta_2 x_1^{-1}x_4,\eta_3 x_1x_2^{-1}x_3^{-1}x_4,\eta_4x_3^{-1}x_4)^t\big) \\ &=  (\eta_2 x_1^{-1}x_4,\eta_3 x_1x_2^{-1}x_3^{-1}x_4+\eta_4x_3^{-1}x_4)^t.
\end{align*}
The supports of the columns of $A_1$ are given by $\alpha_1=\{2\}$, $\alpha_2=\{1,2\}$.
Then, the description of $\mmG_{A,V}$ according to \eqref{eq:AV} is
$$\mmG_{A,V} = \{ f_\eta \in \mmF_{\widehat{A},\widehat{V}} \st \eta_{2} = \eta_{4}\},$$
which agrees with \eqref{eq:geta2}.

\subsection{Gauss reduction and multiple zeros}\label{sec:multi}
There exist families $\mmF_{A,V}$ that do not have multiple $S$-zeros but for which the families $\mmF_{A',V'}$ are not $S$-injective for all possible column permutations of $A$. Therefore, the steps outlined above do not guarantee the existence of multiple $S$-zeros for some member of $\mmF_{A,V}$.
However, under some extra conditions, the existence of multiple $S$-zeros can be asserted. 
These conditions are described in this section.

\paragraph{Example 2. } Let $A$ be as in \eqref{eq:Ared} and let $V$ be
\begin{equation}\label{eq:Vex2}
 V={\small \left(\begin{array}{rccr} 1 & 0 & 1 & -1  \\  1 & 1 & 0 & 0   \\  0 & 1 & 0 & 0  \\   0 & 0 & 1 & 1 \end{array}\right)}.  
\end{equation}
With this matrix $V$, the map $g_{\wk}(x)$ in \eqref{eq:geta2}
 becomes
\begin{equation}\label{eq:g2}
g_{\k_3,\k_4}(x)= \begin{pmatrix}  1 & 0 & 0  \\  0 & 1 & 1 \end{pmatrix}   \begin{pmatrix}  \k_4 x_1^{-1} x_4  \\  \k_3 x_1 x_2^{-1}x_3^{-1} x_4  \\  \k_4  x_1^{-1}x_2^{-1}x_3^{-1}x_4 \end{pmatrix}. 
\end{equation}
Then, $\mmG_{A,V}$ is included in $\mmF_{A',V'}$ with $A'$  and $V'=(v_{ij}')$ given in \eqref{eq:CW}, except for the last column of $V'$ where $v_{13}'=v_{23}'=-1$.
With $S$ given in \eqref{eq:S}, the family $\mmF_{A',V'}$ is not $S$-injective. 
Indeed, for a suitable $Z$, the matrix $M_{\k,\lambda}$ is
$$M_{\k,\lambda} = {\small \begin{pmatrix}  1 & 1 & 0 & 0  \\  0 & 0 & 1 & 1 \\ \k_1\lambda_1 & 0 & 0 & \k_1\lambda_4 \\
(\k_2-\k_3)\lambda_1 & -(\k_2+\k_3)\lambda_2 & -(\k_2+\k_3)\lambda_3 & (\k_2+\k_3)\lambda_4
\end{pmatrix}}. $$
We have that
$$\det(M_{\k,\lambda})=\k_{1}(\k_{2}\lambda_{1}\lambda_{3}+2\,\k_{3}\lambda_{1}\lambda_{4}-\k_{2}\lambda_{2}\lambda_{4}+\k_{3}\lambda_{1}\lambda_{3}-\k_{3}\lambda_{2}\lambda_{4}),
$$
which has both positive and negative coefficients. Therefore, there exist $\eta\in \R^3_{>0}$ and distinct $x,y\in \R^4_{>0}$ such that $x-y\in S$ and  
$ A'\big(\eta \circ x^{V'}\big) = A'\big(\eta \circ y^{V'}\big).$ That is,
\begin{align*}
\eta_1 x_1^{-1}x_4 &=  \eta_1 y_1^{-1}y_4    \\  
 \eta_2 x_1 x_2^{-1}x_3^{-1} x_4 + \eta_3  x_1^{-1}x_2^{-1}x_3^{-1}x_4 &= \eta_2 y_1 y_2^{-1}y_3^{-1} y_4+  \eta_3  y_1^{-1}y_2^{-1}y_3^{-1}y_4.
\end{align*}
The first equality holds independently of $\eta_1$. Hence for $\k_3=\eta_2$ and $\k_4=\eta_3$, we have
$ g_{\k_3,\k_4} (x)= g_{\k_3,\k_4}(y).$
We define $\k_1,\k_2>0$ as
$$ \begin{pmatrix} \k_1 \\ \k_2 \end{pmatrix} = g_{\k_3,\k_4} (x)\quad (= g_{\k_3,\k_4}(y)). $$
Reversing the steps from $f_\k$ to $g_{\wk}$, it follows that
$$A(\k \circ x^V)=A(\k \circ y^V)=0,\qquad x-y\in S,\quad x\neq y, $$
and thus $\mmF_{A,V}$ has multiple $S$-zeros.

\paragraph{Asserting multiple $S$-zeros. } Assume that 
there exist parameters 
$\wk\in \R^{r-s}_{>0}$ and distinct $x,y\in \R^n_{>0}$ with $x-y\in S$ such that 
\begin{align}
g_{\wk} (x) &= g_{\wk}(y),  \label{eq:g} \\
g_{\wk} (x) &\in\R^{s}_{>0}. \label{eq:pos} 
\end{align}
Define $\k\in \R^{r}_{>0}$ by $(\k_1,\dots,\k_s)=g_{\wk}(x)$ and $(\k_{s+1},\dots,\k_r)=\wk$.
Then we have that $A(\kappa\circ x^V)=A(\kappa\circ y^V)=0,$ and hence the family $\mmF_{A,V}$ has multiple $S$-zeros.

To check whether there exist $\wk,x,y$ such that \eqref{eq:g} and \eqref{eq:pos} are fulfilled, we follow the following strategy.

\paragraph{\normalfont \it Checking {\small \eqref{eq:g}}. }
Assume that $\mmF_{A',V'}$ is not $S$-injective.  Then neither is $\mmF_{\widehat{A},\widehat{V}}$ and there exist parameters $\eta\in \R^{s(r-s)}_{>0}$ and distinct $x,y\in \R^n_{>0}$ with $x-y\in S$ such that 
$f_\eta(x) = f_\eta (y)$, for $f_\eta \in \mmF_{\widehat{A},\widehat{V}}$.
We can manually inspect the form of such a vector $\eta$ and obtain one such that $f_\eta$ belongs to $G_{A,V}$.

A sufficient and implementable condition for \eqref{eq:g}  is given by the next proposition. 
The idea is that in order to conclude that \eqref{eq:g} holds for some $\wk\in \R^{r-s}_{>0}$, we need sufficient freedom  to modify   $\eta$ to satisfy the relations  in \eqref{eq:AV}.

By modify we mean the following. Given $\epsilon \in \R^{s}_{>0}$, 
define $\widehat{\epsilon}\in \R^{s(r-s)}_{>0}$ 
by 
$$\widehat{\epsilon}_\ell = \epsilon_{i},\qquad \textrm{if }\quad\ell=j+(i-1)(r-s),\quad \textrm{with } i\in [s],\ j\in [r-s],$$ 
that is,
$ \widehat{\epsilon} =(\epsilon_1,\dots,\epsilon_1, \dots, \epsilon_s,\dots,\epsilon_s)$
such that each $\epsilon_i$ is repeated $r-s$ times.

Let $\widehat{\eta} =\widehat{\epsilon}\circ \eta \in \R^{s(r-s)}_{>0}$.
  If $f_\eta \in \mmF_{\widehat{A},\widehat{V}}$, then $\epsilon f_\eta = f_{\widehat{\eta}}\in \mmF_{\widehat{A},\widehat{V}}$.
Therefore, if $\widehat{\eta}$ fulfils the conditions in  \eqref{eq:AV}, then $ f_{\widehat{\eta}} \in \mathcal{G}_{A,V}$,
 $f_{\widehat{\eta}}(x)=f_{\widehat{\eta}}(y)$, and \eqref{eq:g} is fulfilled.

\begin{proposition} \label{prop:g} Let $\eta\in   \R^{s(r-s)}_{>0}$. 
Assume that for each $j\in [r-s]$ there exists $\ell_j\in [s]$ such that $ \alpha_j \cap \alpha_{j'} =\{\ell_j\} = \{\ell_{j'}\}$ for all $j,j'$ such that the cardinality of $\alpha_j,\alpha_{j'}$ is at least two and such that 
$ \alpha_j \cap \alpha_{j'} \neq \emptyset$.

Then there exists $\epsilon \in \R^{s}_{>0}$ such that  the vector $\widehat{\epsilon} \circ \eta \in \R^{s(r-s)}_{>0}$ 
fulfils the conditions in  \eqref{eq:AV}.
\end{proposition}
\noindent
{\bf Proof. } 
Let $b\in [s]$. If $b\in  \bigcup_{\{j\st \# \alpha_j>1\} }\alpha_j $, choose $j$ such that $b\in \alpha_j$ and define
$$ \epsilon_b = \frac{\eta_{j+(\ell_j-1)(r-s)}}{\eta_{j+(b-1)(r-s)}}.$$
This is well defined because if $b\in  \alpha_j \cap  \alpha_{j'}$ then 
$b=\ell_j=\ell_{j'}$ and $\epsilon_b=1$.
Define $\epsilon_b=1$ otherwise.
Let $\widehat{\eta} =\widehat{\epsilon}\circ \eta$. 
If $b\in \alpha_j $ and $\# \alpha_j>1$, then 
$$\widehat{\eta}_{j+(b-1)(r-s)} =  \frac{\eta_{j+(\ell_j-1)(r-s)}}{\eta_{j+(b-1)(r-s)}}  \eta_{j+(b-1)(r-s)} =
\eta_{j+(\ell_j-1)(r-s)} $$
and hence the conditions in \eqref{eq:AV} are fulfilled.
\qed
 
 \medskip
 Therefore, if $\mmF_{A',V'}$ is not $S$-injective and the assumptions in Proposition~\ref{prop:g} are fulfilled, then  there exists $\wk\in \R^{r-s}_{>0}$ such that \eqref{eq:g} holds. 
The conditions of the statement of Proposition~\ref{prop:g} are clearly fulfilled if the columns of $A_1$ have disjoint supports.

\paragraph{\normalfont \it Checking {\small \eqref{eq:pos}}. }
If the non-zero entries of $A_1$ are all  negative and each row has at least one negative entry, then 
condition \eqref{eq:pos} is fulfilled for all $\wk$. This is the case in the example above. If this is   not   the case, 
then we can resort to the following result, which is adapted from Lemma 4.1 in \cite{craciun-feinbergI}.

\begin{proposition}\label{prop:pos} Let $A'\in \R^{s\times r'}$, $V'\in \R^{n\times r'}$,  $S\subseteq \R^n$ a vector subspace, and $M_{\k,\lambda}$ as defined in \eqref{eq:mkl} for some choice of $Z$.
Let $\eta\in \R^{r'}_{>0}$ such that 
$$A'\eta\in \R^{s}_{>0}\quad\textrm{ and }\quad\det(M_{\eta,\lambda})=0,\quad \textrm{for some }\lambda\in \R^n_{>0}.$$
Then there exist distinct $x,y\in \R^n_{>0}$ and $\k\in \R^{r'}_{>0}$ such that $x-y\in S$ and
$$ A'(\k\circ x^{V'}) = A'(\k\circ y^{V'})  \in \R^s_{>0}.$$
\end{proposition}
\noindent
{\bf Proof. } The matrix $A' \diag(\eta) V'^{t} \diag(\lambda)$  has nontrivial kernel in $S$, since $\det(M_{\eta,\lambda})=0$. That is, there exists $\gamma\in \ker(A' \diag(\eta) V'^{t} \diag(\lambda))\cap S$, $\gamma\neq 0$.  
We define  $x,y\in \R^n_{>0}$ and $\kappa\in \R^{r'}_{>0}$ by 
\begin{align*}
x_i  &=\begin{cases} \gamma_i/(e^{\gamma_i\lambda_i} -1)  & \textrm{ if }   \gamma_i\neq 0  \\ 1 &   \textrm{ otherwise} \end{cases}\\
 y_i & =x_i e^{\gamma_i\lambda_i}\\
  \kappa &=\big(\eta \circ {V'^t} (\lambda\circ \gamma)\big)/\big(y^{V'}-x^{V'}\big),
  \end{align*}
 where division is componentwise  and $0/0=1$. It is easy to check that $y-x=\gamma\in S$,  $\k\in \R^{r'}_{>0}$, and
$ A'\big(\kappa\circ \big(y^{V'}-x^{V'}\big)\big)=0$, cf. \cite[Section~7]{Feinberg1988}, \cite[Th. 5.6]{feliu-inj}.
By replacing $\gamma$ by $\epsilon \gamma$ for $\epsilon>0$ in the definitions above, define analogously $x_\epsilon,y_\epsilon$.
We have
\begin{align*}
A'\big(\kappa\circ x_\epsilon^{V'}\big) &= A' \Big(\frac{\eta \circ V'^{t} (\lambda\circ \epsilon\gamma) \circ x_\epsilon^{V'}}{y_\epsilon^{V'}-x_\epsilon^{V'}} \Big) 
= A'  \Big(\frac{\eta \circ V'^{t} (\lambda\circ \epsilon\gamma)  }{e^{V'^t (\lambda\circ \epsilon\gamma) } -1 } \Big) \\ & \xrightarrow{\epsilon\rightarrow 0} 
A' (\eta \circ (1,\dots,1)^t) = A'\eta \in \R^{s}_{>0}. 
\end{align*}
Therefore, for $\epsilon$ small enough, we get the desired result.
\qed

\subsection{The method}\label{sec:method}
The procedures described in subsections \ref{sec:math}\ref{sec:inj} and \ref{sec:math}\ref{sec:multi}  give a new injectivity-based method to determine whether a set of functions defined by generalised polynomial maps in the positive orthant admit more than one zero on the positive part of $x^*+S$  for varying $x^*$.

Specifically, given $A,V,S$, we proceed as follows:
\begin{itemize}
\item[\bf 0.] Check whether $\mmF_{A,V}$ is $S$-injective using the determinant criterion with the matrix \eqref{eq:mkl}. If yes the family $\mmF_{A,V}$ does not admit multiple $S$-zeros and stop. If not, proceed.
\item[\bf 1.] Compute the Gauss reduction of $A$ and the function $g_{\wk}$.
\item[\bf 2.] Identify matrices $A',V'$ such that $g_{\wk}\in \mmF_{A',V'}$.
\item[\bf 3.] Check whether the family $ \mmF_{A',V'}$ is $S$-injective. If  yes, the family $\mmF_{A,V}$ does not admit multiple $S$-zeros and stop. If not, proceed.
\item[\bf 4.] Check whether the assumptions in Proposition~\ref{prop:g} are fulfilled. If not, go to step {\bf 7}. If yes, proceed.
\item[\bf 5.] If all nonzero entries of $A_1$ in \eqref{eq:B} are negative and each row contains a nonzero entry, then the family $\mmF_{A,V}$ has multiple $S$-zeros and stop.
If not, proceed.
\item[\bf 6.] Check whether the assumptions in Proposition~\ref{prop:pos} are fulfilled. If yes, the family $\mmF_{A,V}$ has multiple $S$-zeros and stop.  If not, proceed.
\item[\bf 7.] Permute the columns of $A,V$ and the entries of $\k$ simultaneously such that the first $s\times s$ minor of $A$ is nonzero and go back to step {\bf 1}.
\end{itemize}

We refer to the step {\bf 0.} as the \emph{standard injectivity method}, since it is the approach underlying the previous methods. 
For step {\bf 3.} we use the determinant criterion if the rank of $A$ and the dimension of $S$ are appropriate. This is the case in the applications in the next section. 

For step {\bf 7}, we fix an order of the set of subsets of $[r]$ with $s$ elements. At the $i$-th iteration of the method, we consider the permutation that sends the $i$-th subset to the front, such that the set of the first $s$ indices and the set of last $r-s$ indices each remain ordered. We check whether the first $s\times s$ minor of $A$ is nonzero. Initial dependencies of the columns of $A$ can be taken into account beforehand, to reduce the number of checks. For example, if two columns of $A$ are linearly dependent, as is often the case for   reaction networks, then subsets  of $[r]$ with $s$ elements containing the two columns are  disregarded.

Except for step {\bf 6}, which is non-trivial, all the other steps are easily implemented using any mathematical software that allows for symbolic computations. 
We have automatised steps {\bf 1}-{\bf 5} and {\bf 7} in Maple and applied the method to a number of situations in the next section. Only when these steps were inconclusive, was it checked whether step {\bf 6} would give a positive answer (see Section \ref{sec:net}).

 The steps of the method are illustrated for one of the small examples in the next section (cf. \eqref{eq:cat2}).

\section{Application to chemical reaction networks}\label{sec:net}
The main application of the method is to determine whether a 
given chemical reaction network admits multiple steady states. 
We follow  the formalism of CRNT \cite{feinbergnotes,gunawardena-notes}.

\paragraph{Setting. }
A \emph{(chemical) reaction network} over a set $\mathcal{X} = \{X_1,\dots,X_n\}$ is a finite collection of reactions
$$ \sum_{\ell=1}^n \mu_{\ell i} X_\ell \rightarrow \sum_{\ell=1}^n \beta_{\ell i} X_\ell,  \qquad i\in [r], $$
 where $\mu_{\ell i},\beta_{\ell i} \in \Z_{\geq 0}$ and the two sides of a reaction are different. Let $A$ be the $n\times r$ matrix whose $(\ell,i)$-th entry is $\beta_{\ell i}-\mu_{\ell i}$, that is,
the net production of  $X_\ell$ in the $i$-th reaction.

The elements $X_\ell$ correspond to chemical species. We denote the concentration of 
$X_\ell$ by $x_\ell$. The vector of concentrations is $x=(x_1,\dots,x_n)\in \R^n_{\geq 0}$ and the concentration at time $t$ is denoted by $x(t)$, although reference to time is often omitted from the notation.

It is custom to model the evolution of the concentrations in time using ODEs.
A typical choice of  ODE system is based on  so-called \emph{mass-action kinetics}, which is a special case of \emph{power-law kinetics}. 
In this setting, the 
 $j$-th reaction is assigned a vector $v^j\in \R^{n}$, and the monomial $x^{v^j}$ is assumed to be proportional to the rate of the $j$-th reaction when the system has concentration $x$.
 
Let $V$ be the $n\times r$ matrix with columns $v^1,\dots,v^r$ and let $\k\in \R^r_{>0}$ be fixed constants. Then 
the system of ODEs is given by 
 \begin{equation}\label{eq:ODE}
\frac{d x}{dt} = f_\k(x),\qquad f_\k(x)= A\big(\k \circ x^V\big).
\end{equation}
This ODE system is defined for $x\in \Omega_V$, where 
 $\R^n_{>0}\subseteq \Omega_V\subseteq \R^n_{\geq 0}$ is obtained by removing from $\R^n_{\geq 0}$ the $\ell$-th hyperplane orthant whenever the $\ell$-th row of $V$ contains a negative entry.
The vector $\k$ is called the vector of reaction rate constants.

In mass-action kinetics,   the matrix $V=(v_{ij})$ is defined by $v_{\ell i} = \mu_{\ell i}$. That is, the $i$-th column of $V$ consists of the coefficients in the left-hand side of the  $i$-th reaction.

The subspace $S=\im(A)\subseteq \R^n$  is called the \emph{stoichiometric subspace}. Because the derivative of $x$ belongs to $S$, the trajectories of \eqref{eq:ODE} are confined to  sets $(x^*+S)\cap \R^n_{\geq 0}$ with $x^*\in \R^n_{\geq 0}$  the initial condition of the system. These sets are called \emph{stoichiometric compatibility classes}. 
We study therefore the dynamics of \eqref{eq:ODE} confined to the stoichiometric compatibility classes and determine the steady states within each   stoichiometric compatibility class. 
In this context, if  the family $\mmF_{A,V}$ has multiple $S$-zeros, then one says that 
 the reaction network is \emph{multistationary} or \emph{admits multiple steady states}.

In the following examples we provide a series of reaction networks and use our method to determine whether multiple steady states can exist in some stoichiometric compatibility class for some choice of reaction rate constants $\k$.

\paragraph{Two-component systems. } Consider the  reaction network
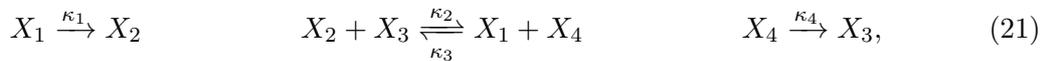
\begin{align}\label{eq:HK}
X_1 & \cee{->[\k_1]} X_2 &  X_2 + X_3  & \cee{<=>[\k_2][\k_3]} X_1 + X_4  &  X_4  & \cee{->[\k_4]}  X_3,
\end{align}
where, as it is custom,  reaction rate constants are written as reaction labels. 
This network models a simple bacterial two-component system in which $X_1,X_2$  (resp. $X_3,X_4$) are the unphosphorylated and phosphorylated forms of a histidine kinase (resp. response regulator).
Using mass-action kinetics, the ODE system \eqref{eq:ODE} of network \eqref{eq:HK} has matrices
$$A =  {\small  \left(\begin{array}{rrrr} -1 & 1 & -1 & 0 \\  1 & -1 & 1 & 0  \\  0 & -1 & 1 & 1 \\  0 & 1 &  -1 & -1  \end{array}\right)},\qquad
 V={\small \left(\begin{array}{cccc} 1 & 0 & 1 & 0  \\  0 & 1 & 0 & 0   \\  0 & 1 & 0 & 0  \\   0 & 0 & 1 & 1 \end{array}\right)}, $$ 
such that
$$f_\k(x)=  A (\k_1 x_1, \k_2 x_2 x_3, \k_3 x_1 x_4, \k_4 x_4)^t. $$
The family $\mmF_{A,V}$ is $S$-injective, and hence there is not a choice of constants for which a stoichiometric compatibility class has multiple steady states. 

The matrix $A$ has rank $2$, and notice that the second and fourth row  are precisely the matrix $A$ in \eqref{eq:AVex1}.
Further, the vector subspace $S$ in \eqref{eq:S} is precisely the stoichiometric subspace of this network, that is, the image of $A$.
The matrix of exponents $V$ in \eqref{eq:AVex1} provides another kinetics for  network \eqref{eq:HK}, which is not mass-action. 
By the results above,  network \eqref{eq:HK}  with the kinetics given by $V$ is not multistationary.

On the other hand, we have shown that with   $V$ given in \eqref{eq:Vex2},  network \eqref{eq:HK} admits
multiple steady states in one stoichiometric compatibility class for some choice of reaction rate constants $\k\in \R^4_{>0}$. Hence the network is multistationary.

 \paragraph{The Langmuir-Hinselwood mechanism. }
 Our second example is  the catalytic oxidation of CO on a Pt(111) surface, which follows the Langmuir-Hinselwood mechanism, and is known to admit multiple steady states. The Langmuir-Hinselwood generally describes the adsorption of one or more reactants on a surface, cf. \cite[Eq. (19)]{razon:1987dv}.
 
Using \cite[Eq. (1)]{Ertl:1980hz} and \cite[Eq. (19)]{razon:1987dv}, the catalytic oxidation of CO is described by the reactions
 \begin{align*}
{\rm CO}+ S & \cee{<=>[\k_1][\k_2]} {\rm CO}_{\rm ad} &  {\rm O}_2 + 2S & \cee{->[\k_3]} 2{\rm O}_{\rm ad} &  {\rm CO}_{\rm ad}+{\rm O}_{\rm ad}   & \cee{->[\k_4]}  {\rm CO}_2+2S,
\end{align*}
 where $S$ represents the active catalytic site on the surface.
 The gases CO, O$_2$ and CO$_2$ are assumed constant. By letting $X_1=S$, $X_2=$CO$_{\rm ad}$, and $X_3=$O$_{\rm ad}$, the reaction scheme  is thus reduced to
  \begin{align}\label{eq:cat2}
X_1 & \cee{<=>[\k_1][\k_2]} X_2&   2X_1 & \cee{->[\k_3]} 2X_3 &  X_2+X_3   & \cee{->[\k_4]}  2X_1.
\end{align}
Assuming mass-action kinetics, the ODE system  in \eqref{eq:ODE} modelling the evolution of the concentrations of $X_1,X_2,X_3$ in time is
\begin{align*}
\frac{dx_1}{dt} &= -\k_1 x_1 + \k_2 x_2 - 2\k_3 x_1^2 + 2\k_4 x_2x_3, \\
\frac{dx_2}{dt}  &= \k_1 x_1 - \k_2 x_2- \k_4 x_2x_3, \\
\frac{dx_3}{dt} &= \k_3 x_1^2 - \k_4 x_2x_3. 
\end{align*}
The matrices $A,V$ in \eqref{eq:ODE}  are
$$A =  {\small  \left(\begin{array}{rrrr} -1 & 1 & -2 & 2\\  1 & -1 & 0 & -1 \\  0 & 0   & 2 & -1  \end{array}\right)},\qquad
 V={\small \left(\begin{array}{cccccc} 1 & 0 & 2 & 0  \\  0 & 1 & 0  & 1   \\  0 & 0 & 0 & 1  \end{array}\right)}. $$ 
We let $S=\im(A)$.
The matrix $A$ has rank $2$. Since the first two rows of $A$ are linearly independent, we redefine $A$ to be the submatrix of $A$ that consists of the first two rows. 
We now go through steps {\bf 0.}-{\bf 7.} of the method. 
The matrix $M_{\kappa,\lambda}$  in \eqref{eq:mkl} is
$$M_{\kappa,\lambda} =  {\small  \left(\begin{array}{ccc} 1 & 1 & 1 
\\  \k_1\lambda_1 & -\k_2\lambda_2 - \k_4\lambda_2 & -\k_4\lambda_3  \\  4\k_3\lambda_1 & -\k_4\lambda_2   &  -\k_4\lambda_3  \end{array}\right)}$$
and its determinant is
$$ -\k_1\k_4\lambda_1\lambda_2+\k_1\k_4\lambda_1\lambda_3+4\k_2\k_3\lambda_1\lambda_2+\k_2\k_4\lambda_2\lambda_3+4\k_3\k_4\lambda_1\lambda_2-4\k_3\k_4\lambda_1\lambda_3.$$
Since the determinant, seen as a polynomial in $\k,\lambda$, has coefficients with opposite signs, the family $\mmF_{A,V}$ is not $S$-injective  (step {\bf 0.}). The Gauss reduction of $A$ is the matrix
$$ {\small  \left(\begin{array}{rrrr} 1 & -1 & 0 & -1 \\  0 & 0 & 1 & -\frac{1}{2}  \end{array}\right)}, \quad\textrm{from where we find} \quad
A_1 = {\small  \left(\begin{array}{rr}  -1  & -1 \\    0  & -\frac{1}{2}  \end{array}\right)},  $$
after permuting the second and third columns of the Gauss reduction of $A$.
The map $g_{\widehat{\k}}(x)$ in \eqref{eq:geta}  becomes (step {\bf 1.}): 
$$ g_{\widehat{\k}}(x)  =  \begin{pmatrix}  \k_2 x_1^{-1}x_2 + \k_4 x_1^{-1}x_2x_3 \\   \frac{1}{2}\k_4 x_1^{-2}x_2x_3 \end{pmatrix}   =  {\small  \left(\begin{array}{rrr}  1 & 1 & 0\\  0 & 0   & \frac{1}{2}  \end{array}\right)}
 \begin{pmatrix}  \k_2 x_1^{-1}x_2 \\ \k_4 x_1^{-1}x_2x_3 \\   \k_4 x_1^{-2}x_2x_3 \end{pmatrix}.$$
The matrices $A',V'$ are thus (step {\bf 2.}):
$$ A'=  {\small  \left(\begin{array}{rrr}  1 & 1 & 0\\  0 & 0   & \frac{1}{2}  \end{array}\right)},
\qquad V'={\small \left(\begin{array}{rrr} -1 & -1 & -2  \\  1 & 1   & 1   \\  0  & 1 & 1  \end{array}\right)}.$$
We compute the new matrix $M_{\kappa,\lambda}$ from $A',V'$ and $S$:
$$M_{\kappa,\lambda} =  {\small  \left(\begin{array}{ccc} 1 & 1 & 1 
\\  -\k_1\lambda_1 + \k_2\lambda_1  & -\k_1\lambda_2 + \k_2\lambda_2 & -2\k_1\lambda_3 + \k_2\lambda_3  \\ 0 & \frac{1}{2}\k_3\lambda_2   &  \frac{1}{2} \k_3\lambda_3  \end{array}\right)}$$
and its determinant 
$$\frac{1}{2} \k_3 ( \k_1 \lambda_2 \lambda_3+ \k_1   \lambda_1 \lambda_3-\k_2   \lambda_1 \lambda_3- \k_1  \lambda_1 \lambda_2+ \k_2   \lambda_1 \lambda_2).$$
Since the determinant has coefficients with opposite sign, the family $\mmF_{A',V'}$ is not $S$-injective (step {\bf 3.}).   Since the support of each of the columns of $\widehat{A}=A'$ contains only one element, the assumptions in Proposition~\ref{prop:g} are fulfilled (step {\bf 4.}). The nonzero entries of $A_1$ are all negative and each row has nonzero entries. By  step {\bf 5.} we conclude that the family $\mmF_{A,V}$ admits multiple $S$-zeros and therefore the network admits multiple steady states in some stoichiometric compatibility class.

\paragraph{Bifunctional kinase. }
Consider the following  reaction network:
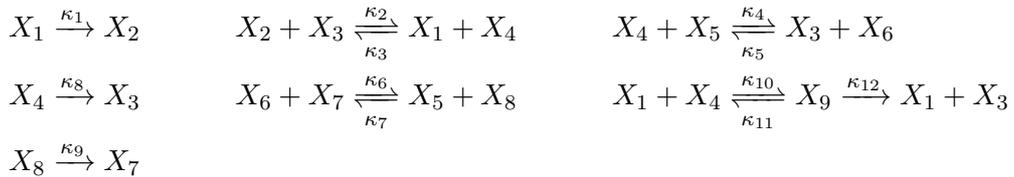
\begin{align*}
X_1 & \cee{->[\k_1]} X_2    & X_2+X_3 & \cee{<=>[\k_2][\k_3]} X_1+X_4  & X_4+X_5 & \cee{<=>[\k_4][\k_5]} X_3+X_6 
\\
 X_4 & \cee{->[\k_8]} X_3 & X_6+X_7 & \cee{<=>[\k_6][\k_7]}  X_5 + X_8   &  X_1 + X_4 &\cee{<=>[\k_{10}][\k_{11}]} X_9\cee{->[\k_{12}]}X_1+X_3
 \\ 
   X_8  &\cee{->[\k_9]} X_7  
   \end{align*}
This network models a phosphorelay signalling system with bifunctional kinase \cite{varun-13}.
Using mass-action kinetics, the matrices $A,V$ in \eqref{eq:ODE} are given as
{\small $$A= \left(  \begin{array}{rrrrrrrrrrrr} -1&1&-1&0&0&0&0&0&0&-1&1&1
\\ 1&-1&1&0&0&0&0&0&0&0&0&0\\ 0&-1
&1&1&-1&0&0&1&0&0&0&1\\ 0&1&-1&-1&1&0&0&-1&0&-1&1&0
\\ 0&0&0&-1&1&1&-1&0&0&0&0&0\\ 0&0
&0&1&-1&-1&1&0&0&0&0&0\\ 0&0&0&0&0&-1&1&0&1&0&0&0
\\ 0&0&0&0&0&1&-1&0&-1&0&0&0\\ 0&0
&0&0&0&0&0&0&0&1&-1&-1\end {array}
 \right)
$$}
and
$$V={\small  \left(\begin {array}{cccccccccccc} 1&0&1&0&0&0&0&0&0&1&0&0
\\ 0&1&0&0&0&0&0&0&0&0&0&0\\ 0&1&0
&0&1&0&0&0&0&0&0&0\\ 0&0&1&1&0&0&0&1&0&1&0&0
\\ 0&0&0&1&0&0&1&0&0&0&0&0\\ 0&0&0
&0&1&1&0&0&0&0&0&0\\ 0&0&0&0&0&1&0&0&0&0&0&0
\\ 0&0&0&0&0&0&1&0&1&0&0&0\\ 0&0&0
&0&0&0&0&0&0&0&1&1\end {array} 
 \right)}.
$$
The family $\mmF_{A,V}$ is not $S$-injective. Since the rank of the stoichiometric matrix $A$ is $5$, we redefine $A$ to be the submatrix of $A$ consisting of rows $2, 4, 6, 8, 9$. Using steps {\bf 1}-{\bf 3}, we  conclude that the network does not admit multiple steady states.
In other words, Gauss reduction of this submatrix of $A$ gives two matrices $A'$, $V'$ such that the family $\mmF_{A',V'}$  is $S$-injective.

\paragraph{Apoptosis. }
We next consider a reaction network, which is a basic model of caspase activation for apoptosis \cite{Eissing2004}:
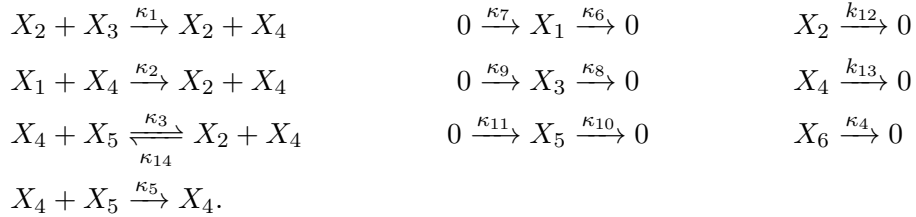
\begin{align*}
X_2+X_3 & \cee{->[\k_1]}  X_2+ X_4   &  0 \cee{->[\k_7]} X_1 & \cee{->[\k_6]} 0 &   X_2 & \cee{->[k_{12}]} 0 \\
X_1+X_4 & \cee{->[\k_2]}  X_2+ X_4   &  0  \cee{->[\k_9]} X_3 & \cee{->[\k_{8}]} 0 & X_4&  \cee{->[k_{13}]} 0 \\
X_4+X_5 & \cee{<=>[\k_3][\k_{14}]}  X_2+ X_4   &  0  \cee{->[\k_{11}]} X_5 &\cee{->[\k_{10}]} 0 & X_6 & \cee{->[\k_4]} 0 \\
X_4 + X_5 & \cee{->[\k_5]} X_4. 
  \end{align*}
With mass-action kinetics, this network is known to admit multiple steady states  \cite{Eissing2004}.
The rank of the stoichiometric subspace is maximal, that is, $S=\R^6$.
For this network, steps {\bf 1}-{\bf 5} are inconclusive, for all possible column permutations, that is, we can neither assert nor reject that multiple steady states exist. 

We permute the columns of $A,V$ such that the columns 7, 9, 11, 12, 13, 14 are first 
and apply Gauss reduction. The matrix $A_1$ in \eqref{eq:B} becomes 
$$ A_1={\small \left(\begin {array}{rrrrrrrr}  
0&-1&0&0&0&-1&0&0 \\ 
-1&0&0&0&0&0&-1&0
\\ 0&0&0&-1&-1&0&0&-1
\\ 0&-1&0&0&0&0&0&0\\ 
-1&0&0&1&0&0&0&0\\ 0&0&-1&1&0
&0&0&0\end {array}\right)}.
$$
Observe that $A_1$ fulfils the assumptions in Proposition~\ref{prop:g}. 
We check step {\bf 6}. The matrices $A',V'$ are in this case
$$A'={\small \left( \begin {array}{rrrrrrrrrrrr} 1&1&0&0&0&0&0&0&0&0&0&0
\\ 0&0&1&1&0&0&0&0&0&0&0&0\\ 0&0&0
&0&1&1&1&0&0&0&0&0\\ 0&0&0&0&0&0&0&1&0&0&0&0
\\ 0&0&0&0&0&0&0&0&1&-1&0&0\\ 0&0&0
&0&0&0&0&0&0&0&1&-1\end {array} \right)}
$$
and 
$$V'= {\small \left( \begin {array}{rrrrrrrrrrrr} 1&1&0&0&0&0&0&1&0&0&0&0
\\ 0&0&1&0&0&0&0&-1&1&0&0&0\\ 0&0&
1&1&0&0&0&0&1&0&0&0\\ 1&0&0&0&0&1&0&1&-1&-1&1&0
\\ 0&0&0&0&0&1&1&0&0&0&1&0\\ 0&0&0
&0&1&0&0&0&0&1&-1&0\end {array}  \right)}.
$$
The vector
$\eta = (1,3,1,1,1,1,1,1,16/15,1,1,1/2) $
fulfils $A'\eta\in \R^6_{>0}$ and that the determinant of $M_{\eta,\lambda}$
vanishes for all $\lambda\in \R^6_{>0}$.
Therefore, by Proposition~\ref{prop:pos} and step  {\bf 6}, we conclude that  this network admits multiple steady states.

\paragraph{Networks of gene regulation. } In \cite{siegal-gaskins} a total of 40,680 reaction networks modelling gene regulatory systems  with mass-action kinetics are considered and analysed for multistationarity. The authors use the CRNT toolbox \cite{crnttoolbox}
together with a method termed \emph{network ancestry}, justified by theoretical results in \cite{joshi-shiu-II}.

The authors determine that 2,654 out of the 40,680 networks cannot have multiple steady states. The standard injectivity method correctly precludes multistationarity for 691 of these 2,654 and is inconclusive for the remaining networks. Using  the method presented here we can conclude that  all  2,654 networks are not multistationary.

Interestingly, the authors cannot decide whether multiple steady states occur for 1,050 out of the 40,680 networks.
The method successfully classifies these 1,050 networks and we can conclude that 47 of them are multistationary, while the remaining 1,003 cannot have multiple steady states.
The remaining 36,976 networks  are shown to be multistationary in \cite{siegal-gaskins}. The method is applied  to the smallest 2,000 of these, and we reach the same conclusion.

\paragraph{Atoms of multistationarity. } In \cite{joshi-shiu-II},  all possible networks consisting of two reactions and at most two molecules at each side of a reaction are considered. Here each reaction can be either irreversible $\cee{->}$ or reversible $\cee{<=>}$. So-called flow reactions $0\cee{<=>}  X$ are added for each chemical species, such that $S=\R^n$. The authors  determine which of these networks are multistationary when taken with mass-action kinetics.

We  considered the 142 networks for which the standard injectivity criterion does not rule out multistationarity. Among these, the authors show that precisely 35 are multistationary and hence 107 are not multistationary.
The new method does not perform as good as in the other test cases. Indeed, only 24 networks are identified as non-multistationary, while the remaining 118 are left unclassified.

A plausible explanation for the failure is the following. 
With mass-action kinetics, a reaction of the form $0\rightarrow X$, called inflow, contributes a constant term $\k_i$ to the $i$-th component of $f_\k$ in \eqref{eq:ODE}. Therefore, 
we can write $f_\k(x)= \widetilde{f}_\eta(x) + \widetilde{\k}$, for some vector $\widetilde{\k}\in \R^n_{>0}$. It follows that $f_\k(x)=f_\k(y)$ if and only if $\widetilde{f}_\eta(x)=\widetilde{f}_\eta(y)$ and hence the standard injectivity method checks for injectivity of $\widetilde{f}_\eta$. In fact, $\widetilde{f}_\eta$ is  $g_{\wk}$ in our method, if the reactions are ordered such that the inflow reactions are first.

\section{Discussion}
Injectivity-based methods are used to preclude multistationarity in reaction networks \cite{craciun-feinbergI,craciun-feinberg-semiopen,feliu-inj,banaji-craciun2,muller,soule,PerezMillan,wiuf-feliu,MullerSigns,joshi-shiu-I,shinar-conc}.  The rationale behind the methods  is that multiple zeros cannot occur if the map is injective. 
In \cite{craciun-feinbergI,craciun-feinberg-semiopen,feliu-inj,wiuf-feliu,joshi-shiu-I} the modelling  framework is either mass-action or power-law kinetics, similar to what is used here. A  common aspect of these works is that the mathematical development focuses on the vector subspace $S$ being the image of the coefficient matrix $A$.

In \cite{PerezMillan,muller}, injectivity of a monomial map is studied in order to assert or preclude multistationarity. In both works, the authors are interested in determining whether a generalised binomial map, that is, a generalised polynomial map with two terms, admits multiple positive zeros.  The positive zeros of a binomial map can be parameterised by a monomial map obtained by dividing one term of the binomial by the other term. 
In this way, there is passage from the non-existence of multiple positive zeros of the binomial map to the injectivity of the monomial map, which is identical to the passage from the study of $f_\kappa$ to the study of $g_{\wk}$ here.

The connection between  \cite{PerezMillan,muller} and \cite{craciun-feinbergI,craciun-feinberg-semiopen,feliu-inj,wiuf-feliu} is clarified in \cite{MullerSigns}, where unifying sign and determinant conditions for the injectivity of generalised polynomial maps are given. An important novelty of \cite{MullerSigns} is that $S$ is given independently of the image of $A$, in contrast to earlier work. 
This uncoupling is key in the results presented here.

We have demonstrated by examples that the new method can be applied to a vast amount of networks for which standard injectivity approaches are inconclusive. We have further applied the method to the five basic building blocks   in cell signalling in  \cite{Feliu:royal} that are shown to admit multiple steady states. The method correctly classifies the networks as multistationary as well.
 However, we are not guaranteed that the method will classify any given network, as we discussed for the two-reaction networks in \cite{joshi-shiu-II}.

We  focus on Gauss reduction of the generalised polynomial maps. Gauss reduction preserves  the number of equations, which implies that the determinant criterion can be used to check step {\bf 3}, when applying the method to reaction networks. Further, the method can be applied as a black box. However, the steps presented here also
apply if we replace $f_\k$ by any set of generalised polynomial maps with the same positive zeros as $f_\k$. For example, if $f_\k$ is polynomial, one might consider a Gr\"obner basis of the system. In this sense, our work is a generalisation of \cite{PerezMillan} to the case where the Gr\"obner basis is not binomial.  If the new set of equations differs from $f_\k$ in  number of equations, then the sign criterion is to be used at step {\bf 3}, which can also be computationally checked  \cite{MullerSigns}. 

Finally, the method has been presented in connection with reaction networks. However, the mathematical framework is given in full generality and can be used to obtain information on the number of positive zeros of generalised polynomial maps, independently of the context in which the question arises.
In particular,  since any polynomial can be embedded into a family $\mmF_{A,V}$, the method can be applied to preclude the existence of multiple positive solutions to any polynomial equation by letting $S=\R^n$.

\section*{Acknowledgments}
B. Joshi, A. Shiu,  and D. Siegal-Gaskins are  thanked for providing raw data and discussions of their results in \cite{joshi-shiu-II} and \cite{siegal-gaskins}. C. Wiuf is thanked for useful discussions and comments on earlier versions of this manuscript. My co-authors in \cite{MullerSigns} are thanked for interesting discussions on injectivity.
An anonymous referee is thanked for interesting comments and pointing the author to relevant literature and examples. 

This work has been partially supported by project  MTM2012-38122-C03-01/FEDER from the Ministerio de Econom\'{\i}a y Competitividad, Spain, the Carlsberg Foundation, and the Lundbeck Foundation.


{\footnotesize

\end{document}